\documentclass{article}
\usepackage{latexsym}
\usepackage{amssymb}
\usepackage{amsmath}
\usepackage{layout}
\newtheorem{teo}{Theorem}[section]
\newtheorem{lema}[teo]{Lemma}
\newtheorem{cor}[teo]{Corollary}
\newtheorem{prop}[teo]{Proposition}
\newtheorem{obs}[teo]{Remark}
\newtheorem{defin}{Definition}
\newenvironment{dem}{\text\bf Proof:}{}
\newcommand{\R}{\mathbb{R}}
\newcommand{\N}{\mathbb{N}}

\newcommand{\D}{\mathbb D}

\newcommand{\cf}{\mathcal{F}}
\newcommand{\cs}{{\cal S}}
\newcommand{\cc}{{\cal C}}

\begin{document}

\begin{center}
\huge {

Integration with respect to local time and It{\^o}'s formula for
smooth nondegenerate martingales}

\vspace{.5cm}

\normalsize {\bf Xavier Bardina$^{1,*}$} and {\bf Carles
Rovira$^2$}

{\footnotesize \it $^1$ Departament de Matem\`atiques, Universitat
Aut\`onoma de Barcelona, 08193-Bellaterra (Barcelona), Spain.

$^2$ Facultat de Matem\`atiques, Universitat de Barcelona, Gran
Via 585, 08007-Barcelona, Spain.

 {\it E-mail addresses}: Xavier.Bardina@uab.cat,
Carles.Rovira@ub.edu}

{$^{*}$corresponding author}

\end{center}

\begin{abstract}%
We show an It\^ o's formula for nondegenerate Brownian martingales
$X_t=\int_0^t u_s dW_s$ and functions $F(x,t)$ with locally
integrable derivatives in $t$ and $x$. We prove that one can express
the additional term in It\^o's s formula as an integral over space
and time with respect to local time.

\end{abstract}



{\bf Running title:} It{\^o}'s formula for nondegenerate
martingales

\section*{Introduction}

We consider a continuous nondegenerate martingale $X=\{X_t, \,t \in
[0,1]\}$ of the form $X_t = \int_0^t u_s dW_s$ where $W=\{W_t, t \in
[0,1]\}$ is a standard Brownian motion and $u$ is an adapted
stochastic process. Let $F: \R \times [0,1] \rightarrow \R$ be an
absolutely continuous function with partial derivatives satisfying
some local integrability properties. The main aim of this paper is
to obtain an It\^o's formula for $F(X_t,t)$ where the term
corresponding usually to the second order derivative is expressed as
an integral over space and time with respect to local time.

We will prove this results when $u$ satisfies (locally) the
assumptions

\begin{itemize}

\item[{\bf (H1)}] For all $t \in [0,1]$, $u_t$ belongs to the space $\D^{3,2}$ and for all $p \ge 2$
$$E \vert u_t \vert^p + E \vert D_s u_t \vert^p + E \left( \int_{r \vee s}^1 \vert D_rD_s u_\theta \vert^2 d\theta \right)^{p/2} + E \left( \int_{r \vee s \vee v }^1 \vert D_v D_rD_s u_\theta \vert^2 d\theta \right)^{p/2} \le K_p,$$

\item[{\bf (H2)}] $\vert u_t \vert \ge \rho >0$ for some constant
 $\rho$ and for all $t \in [0,1].$

\end{itemize}

Moret and Nualart (2000) consider an It\^o's formula for this
class of nondegenerate martingales. Their main result reads as follows:

\begin{teo}[Moret and Nualart, 2000]\label{teoMN}
Let $u$ be a process satisfying {\bf (H1)} and  {\bf (H2)}. Set
$X=\int_0^t us dW_s$.
Then for any
 funcion $f \in L^2_{\rm loc} (\R)$ the quadratic covariation $[f(X),X]$ exists and the following It\^{o}'s formula holds
$$
F(X_t)=F(0)+ \int_0^t f(X_s) dX_s + \frac12 [f(X),X]_t,$$ for all
$t \in [0,1]$, where $F(x)=F(0)+\int_0^x f(y) dy.$
\end{teo}
Moret (1999), gave an extension of this last result for functions
$F$ depending also on  $t$. They consider a new hypothesis on
functions $f$:
\begin{itemize}

\item[{\bf (C)}] $f(\cdot,t) \in
L^2_{loc}(\R)$ and for all compact set $K\subset\R$ $f(x,t)$ is
continuous in $t$ as a function of $[0,T]$ to $L^2(K)$

\end{itemize}
Then, their result is the following:

\begin{teo}[Moret, 1999]\label{teoM}
Let $u$ be a process satisfying {\bf (H1)} and  {\bf (H2)}. Set
$X=\int_0^t us dW_s$.
Let $F(x,t)$ be an absolutely continuous function in $x$ such that
the partial derivative $f(\cdot,t)$ satisfies {\bf (C)}. Then, the
quadratic covariation $\left[f(X,\cdot),X\right]$ exists  and the
following It\^{o}'s formula holds
$$F(X_t,t)=F(0,0)+\int_0^tf(X_s,s)dX_s+\frac12\left[f(X,\cdot),X\right]+\int_0^tF(X_s,ds),$$
where
$$\int_0^tF(X_s,ds)\equiv \lim_{n\to+\infty}\sum_{t_i\in D_n, t_i\leq
t} \left(F(X_{t_{i+1}},t_{i+1})-F(X_{t_{i+1}},t_i)\right),$$ exists
uniformly in probability for $(D_n)_n$ a sequence of smooth
partitions of $[0,1]$.
\end{teo}

In these two results, following the ideas of F{\"o}llmer, Protter
and Shiryayev (1995) for the Brownian motion,  the additional term
is written as a quadratic covariation. Bardina and Jolis (1997,
2002) extended the results of F{\"o}llmer {\it et al.} (1995) to the
case of the elliptic and hypoelliptic diffusions.

Nevertheless, it is important to point out the differences between the work of Moret and Nualart (2000) and
F{\"o}llmer {\it et al.} (1995). One of the keys of their proofs is to obtain some a priori estimates on the Riemann sums. In F{\"o}llmer {\it et al.} (1995) these estimates are obtained using the semimartingale expression of the time-reversed Brownian motion and well-known bounds for the density of the Brownian motion. Moret and Nualart (2000)
used another approach, using Malliavin calculus in order to
obtain sharp estimates for the density of the process $X_t$ and avoiding the
time-reversed arguments.

We want to express the quadratic variation term as an
integral with respect to the local time. There are several papers
where the integrals with respect to local time are used in It\^o's
formula. In 1981, Bouleau and Yor obtained the following extension
of the It{\^o}'s formula :
\begin{teo}[Bouleau and Yor, 1981]\label{teoBY}
Let $X=\left(X_t\right)_{t\geq0}$ be a continuous semimartingale
and let $F:\R\longrightarrow\R$ be an absolutely continuous
function with derivative $f$. Assume that $f$ is a mesurable
locally bounded function.
Then:
$$
F(X_t)=F(X_0)+\int_0^tf(X_s)dX_s-\frac12\int_{\R}f(x)d_xL_t^x
$$
where $d_xL_t^x$ is an integral with respect to $x\longrightarrow
L_t^x$.
\end{teo}

Eisenbaum (2000,
2001) defined an integral in time and
space with respect to the local time of the Brownian motion. Using
this  integral, the quadratic covariation in the formula given in
F{\"o}llmer {\it et al.} can be expressed as an integral with
respect to the local time. She obtained the following result:

\begin{teo}[Eisenbaum, 2000 and 2001]\label{teoE}
Let $W=\left(W_t\right)_{0\leq t\leq1}$ be a standard Brownian
motion and $F$ a function defined on $\R\times[0,1]$ such that
there exist first order Radon-Nikodym derivatives $\frac{\partial
F}{\partial t}$ and $\frac{\partial F}{\partial x}$ such that for
every  $A\in\R_{+}$,
$$\int_0^1\int_{-A}^A\left|\frac{\partial F}{\partial t}(x.s)\right|\frac1{\sqrt{s}}dxds<+\infty$$
and
$$\int_0^1\int_{-A}^A \left(\frac{\partial F}{\partial x}(x,s)\right)^2\frac1{\sqrt{s}}dxds<+\infty.$$
Then,
\begin{eqnarray*}
F(W_t,t)=F(W_0,0)+\int_0^t\frac{\partial F}{\partial
x}(W_s,s)dW_s+\int_0^t\frac{\partial F}{\partial
t}(W_s,s)ds-\frac12\int_0^t\int_{\R}\frac{\partial F}{\partial
x}(x,s)dL_s^x.
\end{eqnarray*}
\end{teo}
This result has been extended by Bardina and Rovira (2007) for
elliptic diffusion processes.

In our papers we will follow the ideas Eisenbaum (2000,2001), assuming on the function $F$ the hypothesis considered in Theorem \ref{teoE}. In the papers of Eisenbaum (2000,2001), as well as in F{\"o}llmer {\it et al.}
(1995) or in the extension of Bardina and Rovira (2007), one of
the main ingredients is the study of the time reversed process and the relationship between the quadratic covariation and the forward and backward stochastic integrals. We show that we can adapt the methods of
Eisebaum  without using the
time reversed process and the backward integral. We will follow the methods  of Moret and Nualart
(2000) and we will use Malliavin calculus to obtain the necessary estimates for the Riemann sums .

In our paper, the existence of the quadratic covariation is not
one of our main objectives. Nevertheless, it will be an important  tool in our
computations. We recall its definition.
\begin{defin}
Given two stochastic processes $Y=\{Y_t, t \in [0,1]\}$ and
$Z=\{Z_t, t \in [0,1]\}$ we define their quadratic covariation as
the stochastic process $[Y,Z]$ given by the following limit in
probability, if it exists,
$$
[Y,Z]_t = \lim_n \sum_{ t_i \in D_n, t_i < t } (Y_{t_{i+1}}-
Y_{t_i}) (Z_{t_{i+1}}- Z_{t_i}).$$ where  $D_n$ is a sequence of
partitions of $[0,1]$.
\end{defin}

We will assume that the partitions $D_n$ satisfy
\begin{itemize}
\item[{\bf (M)}] $\lim_n \sup_{t_i \in D_n} (t_{i+1} - t_i )=0, \qquad M:=\sup_n \sup_{t_i \in D_n}
\frac{t_{i+1}}{t_i} < \infty.$
\end{itemize}
We impose this condition in order to avoid certain possibly
exploding Riemann sums.

Other extensions for It{\^o}'s formula has been obtained recently.
Among others, there is the paper  of Dupoiron {\it et al.} (2004)
for uniformly elliptic diffusions and Dirichet processes, the work
of Ghomrasni and Peskir (2006) for  continuous semimartingales, the
paper of  Flandoli, Russo and Wolf (2004) for a Lyons-Zheng process
or the work of Di Nunno, Meyer-Brandis, {\O}ksendal and Proske
(2005) for L\'evy processes.

The paper is organized as follows. In Section 1 we give some basic
definitions and results on Malliavin calculus, recalling some
results obtained in Moret and Nualart (2000). In Section 2 we
define the space where we are able to construct an integral in the
plane with respect to the local time of a nondegenerate Brownian
martingale. Finally, Section 3 is devoted to present our main result the extension
of It{\^o}'s formula.

Along the paper we will denote all the constants by $C, C_p, K$ or
$K_p$, unless they may change from line to line.

\section{Preliminaries}

 Let  $(\Omega,\cf,P)$ be the
canonical probability space of a standard Brownian motion
 $W=\{ W_t, 0 \le t \le 1 \}$, that is,
$\Omega$ is the space of all continuous functions $\omega:[0,1]
\to \R$ vanishing at 0, $P$ is the standard Wiener measure on
$\Omega$ and $\cf$ is the completion of the Borel $\sigma$-field
of $\Omega$ with respect to $P$. Let $H=L^2([0,1])$.

Let $\cs$ be the set of smooth random variables of the form
\begin{equation}\label{A1}
F=f(W_{t_1},\ldots,W_{t_n}),
\end{equation}
$f\in \cc_b^\infty(\R^n)$ and $t_1, \ldots, t_n \in [0,1].$
The Malliavin derivative of a smooth random variable $F$ of the form (\ref{A1}) is the
 stochastic process $\{ D_tF, t \in [0,T]\}$
given by
$$D_t F = \sum_{i=1}^n \frac{\partial f}{\partial x_i} (W_{t_1},\ldots,W_{t_n}) I_{[0,t_i]}(t),\qquad t \in [0,1].$$

The Malliavin derivative of order $N \ge 2$ is defined by iteration, as follows. For $F \in \cs,\, t_1,\ldots,t_N
\in [0,T],$
$$
D^N_{t_1,\ldots,t_N} F = D_{t_1}D_{t_2} \ldots D_{t_N} F.$$
For any real number $p \ge 1$ and any integer $N\ge 1$ we denote by
$\D^{N,p}$ the completion of the set $\cs$ with respect to the norm
$$
\Vert F \Vert_{N,p} = \big[ E( \vert F \vert^p )
+ \sum_{i=1}^N E ( \Vert D^i F \Vert^p_{L^2([0,T]^i)} ) \big]^\frac1p.$$
The domain of the derivative operator $D$ is the space $\D^{1,2}.$

The divergence operator $\delta $ is the adjoint of the derivative operator.
The domain of the operator $\delta $,
denoted by Dom $\delta $, is the set of processes $u \in L^{2}([0,T] \times \Omega)$ such that there exists a square integrable random variable $\delta(u)$ verifying
$$
E(F\delta (u))=E\left( \int_0^1 D_tF u_t dt \right),
$$
for any $F \in \cs$.
The operator $\delta$ is an extension of It\^o's stochastic integral and we will make use of the notation $\delta(u)=\int_0^1 u_s dW_s.$

We will recall some useful results from Moret and Nualart (2000). We
refer the reader to this paper for their proof and a detailed
account of these results. We also refer to Nualart (1995, 2006) for
any other property about operators $D$ and $\delta$.

\begin{prop}\label{2.6MN}
Let $Y$ be a random variable in the space $\D^{1,2}$ such that $\int_a^b ( D_s Y)^2 ds > 0$ a.s. for some $0 \le a < b \le 1.$ Assume that $(DY/\int_a^b (D_sY)^2 ds) I_{[a,b]}$ belongs to Dom $\delta$. Then $Y$  has an absolutely continuous distribution with density $p$ that satisfies the inequality
$$
p(x) \le E \Big\vert \int_a^b \left( \frac{D_sY}{\int_a^b (D_s Y)^2
ds } \right) dW_s \Big\vert.$$
\end{prop}
\begin{dem}
It follows from Proposition 1 and (2.6) in Moret and Nualart (2000).
\hfill $\Box$
\end{dem}

The following Proposition is also a  slight modification of Corollary 2 of Moret-Nualart (2000).

\begin{prop}\label{corol2MN}
Let $Y$ be a random variable in the space $\D^{1,2}$ such that $\int_0^1 ( D_s Y)^2 ds > 0$ a.s. Let $Z$ be a positive square integrable random variable such that $( Z DY / \int_0^1 (Ds Y )^2 ds ) I_{[0,1]}$ belongs to Dom $\delta$. Then, for any $f \in L^2(\R)$, we have
$$
\vert E (f(Y)^2 Z) \vert = \Vert f  \Vert_2^2 E \Big\vert \delta \left( \frac{Z DY}{\Vert DY \Vert_H^2} \right) \Big\vert.$$
\end{prop}
\begin{dem}
See Corollary 2 in Moret-Nualart (2000). The same proof works using a dominated convergence argument.
\hfill $\Box$
\end{dem}

\begin{lema}\label{lema10MN}
Fix $p \ge 1$. Suppose that $u$ satisfies hypotheses {\bf (H1)} and {\bf (H2)}. Let $Z \in \D^{1,2p}$.
Then, we have, for $0 \le a < b \le 1$:
\begin{eqnarray*}
&& E \Big\vert \int_a^b Z \frac{ D_t X_b}{ \int_a^b (D_t X_b )^2 dt} dW_t \Big\vert^p \\
&& \qquad \le C_0 (b-a)^{-p/2} \left( \left( E \vert Z \vert^{2p} \right)^{1/2} +
\left( E \big\vert \int_a^b (D_t Z)^2 dt  \big\vert^{p} \right)^{1/2} \right),
\end{eqnarray*}
where $C_0$ is a constant does not depend on $Z$.
\end{lema}
\begin{dem}
See Lemma 10 in Moret-Nualart (2000).
\hfill $\Box$
\end{dem}

\section{Stochastic integration with respect to local time of the martingale}

Following the ideas of Eisenbaum (2000), we consider first the
space of functions for whose elements we can define a stochastic
integration with respect to local time.

Let $f$ be a measurable function from $\R\times[0,1]$ into $\R$.
We define the norm $\| \cdot \|$ by
$$
\|f\|=\left(\int_0^1\int_{\R}
f^2(x,s)\frac1{s^{\frac34}}dxds\right)^{\frac12}
$$

\medskip

Consider the set of functions $$\mathcal
H=\{f:\,\|f\|<+\infty\}.$$ It is easy to check that $\mathcal
H$  is a Banach space.


Let us consider $X$ a nondegenerate martingale of the type
$X_t=\int_0^t u_s dW_s$ where $u$ is an adapted stochastic process
satisfying hypotheses {\bf (H1)} and {\bf (H2)}. Let us show now how
to define a stochastic integration over the plane with respect to
the local time $L$ of the process $X$ for the elements of $\mathcal
H$.

Let $f_{\Delta}$ be an elementary function,
$$f_{\Delta}(x,s):=\sum_{(x_k,s_l)\in\Delta}
f_{kl}I_{(x_k,x_{k+1}]}(x)I_{(s_l,s_{l+1}]}(s),$$where
$(x_k)_{1\leq k\leq m_1}$ is a finite sequence of real numbers,
$(s_l)_{1\leq l\leq m_2}$ is a subdivision of $[0,1]$, \break
$(f_{kl})_{1\leq k\leq m_1;\, 1\leq l\leq m_2}$ is a sequence of
real numbers and finally, $\Delta=\{(x_k,s_l),\,1\leq k\leq
m_1,\,1\leq l\leq m_2\}$. It is easy to check that the elementary
functions are dense in $\mathcal H$.

\medskip

We define the integration for the elementary function $f_{\Delta}$
with respect to the local time $L$ of the martingale $X$ as
follows
$$\int_0^1\int_{\R}f_{\Delta}(x,s)dL_s^x=
\sum_{(x_k,s_l)\in\Delta}
f_{kl}(L_{s_{l+1}}^{x_{k+1}}-L_{s_{l}}^{x_{k+1}}-L_{s_{l+1}}^{x_{k}}+L_{s_{l}}^{x_{k}}).$$
Let $f$ be a function of $\mathcal H$. Let us consider
$(f_n)_{n\in\N}$ a sequence of elementary functions converging to
$f$ in $\mathcal H$. We will check that the sequence
$\left(\int_0^1\int_{\R}f_{n}(x,s)dL_s^x\right)_{n\in\N}$
converges in $L^1$ and that the limit does not depend on the
choice of the sequence $(f_n)_{n\in\N}$. So, we will use this
limit as the definition of the integral
$\int_0^1\int_{\R}f(x,s)dL_s^x$.

\medskip

First of all, let us see a previous lemmas.

\begin{lema}\label{lema1}
For any locally bounded Borel measurable function $f$ and any
$t\in(0,1]$ we have
$$\int_{\R}f(a)d_aL^a_t=-\left[f(X),X\right]_t,$$
where $d_aL_t^a$ denotes the integral with respect to
$a\longrightarrow L_t^a$.
\end{lema}
\begin{dem}
It follows easily from  Theorem \ref{teoBY} and Theorem
\ref{teoMN}.
\hfill $\Box$
\end{dem}

\begin{lema}\label{lemafita}
Consider $f_1(x):=I_{(a,b]}(x)$ and $f_2(x):=I_{(c,d]}(x)$, where
$a<b$ and $c<d$ are real numbers.
Then for all $t_i<t_j\leq t$,
\begin{eqnarray*}
&&E\left[f_1(X_{t_{i+1}})f_2(X_{t_{j+1}})(X_{t_{i+1}}-X_{t_{i}})(X_{t_{j+1}}-X_{t_{j}})\right]\\
&\leq&E\left[f_1(X_{t_{i+1}})f_2(X_{t_{j+1}})C_{ij}\right],
\end{eqnarray*}
where
$$\|C_{ij}\|_2\leq C\frac{(t_{i+1}-t_i)(t_{j+1}-t_j)}{\sqrt{t_{i+1}(t_{j+1}-t_{i+1})}},$$
and $C$ does not depend on $f_1$ and $f_2$.
\end{lema}
\begin{dem} When $f_1=f_2=f\in\mathcal C_K^{\infty}(\R)$, this inequality is checked in the proof of Proposition 14 in Moret and Nualart (2000). The same proof also works when $f_1\neq f_2$ with $f_1,f_2\in\mathcal C_K^{\infty}(\R)$.
Now, fixed our functions $f_1,f_2$ let us consider sequences $f_1^n\uparrow f_1$
and $f_2^n\uparrow f_2$ with $f_n^i\in\mathcal C_K^{\infty}(\R)$ for all
$n$ and $i\in\{1,2\}$. Then, the result can be obtained by a dominated convergence argument.

\hfill $\Box$
\end{dem}

\smallskip

\begin{teo}\label{teo2}
Let $f$ be a function of $\mathcal H$. Then, there exists the
integral  $\int_0^t\int_{\R}f(x,s)dL_s^x$
for any  $t\in[0,1]$.
\end{teo}
\begin{dem}
Let $f_{\Delta}$ be an elementary function. From Theorem \ref{teoM} and Lemma \ref{lema1}
it is easy to get that the quadratic covariation $[f(X,.),X]_t$ exists and that
$$\int_0^t\int_{\R}f_{\Delta}(x,s)dL_s^x=-\left[f_\Delta(X,\cdot),X\right]_t.$$

The key of the proof is to check that for all elementary function $f_\Delta$
\begin{equation}\label{simpli}
E\left(\left|\int_0^t\int_{\R}f_{\Delta}(x,s)dL_s^x\right|\right)\leq C\|f_{\Delta}\|,
\end{equation}
where the constant does not depend on $f_\Delta$.

Notice that,
\begin{eqnarray}
&&E\left(\left|\int_0^t\int_{\R}f_{\Delta}(x,s)dL_s^x\right|\right)\nonumber\\
&=&E\left(\left|\left[f_\Delta(X,\cdot),X\right]_t\right|\right)\nonumber\\
&=& E \left(\left|\lim_{n\to\infty}\sum_{t_i\in D_n,t_i\leq
t}\left(f_{\Delta}(X_{t_{i+1}},t_{i+1})-f_{\Delta}(X_{t_{i}},t_{i})\right)\left(X_{t_{i+1}}-X_{t_i}
\right) \right| \right) \nonumber\\&\leq& \left\{ 2
\liminf_{n\to\infty} E \left( \left|\sum_{t_i\in D_n,t_i\leq
t}f_{\Delta}(X_{t_{i+1}},t_{i+1})\left(X_{t_{i+1}}-X_{t_i} \right)
\right|^2 \right) \right.\nonumber\\&& +\left. 2
\liminf_{n\to\infty} E \left(
\left|\sum_{t_i\in D_n,t_i\leq
t}f_{\Delta}(X_{t_{i}},t_{i})\left(X_{t_{i+1}}-X_{t_i}
\right) \right|^2 \right)  \right\}^{\frac12} \nonumber\\
&:=&
\left( 2 \left(\liminf_{n\to\infty}I_1\right)+2 \left(\liminf_{n\to\infty}I_2\right)\right)^{\frac12},\label{aaa1}
\end{eqnarray}
where in the last inequality we have used Fatou's lemma.

Along the study of $I_1$ and $I_2$ we will make use of the methods
presented in the proofs of Propositions 13 and 14 in Moret and
Nualart (2000). For the sake of completeness, we will give the
main steps of our proofs in the study of $I_2$. For the other
terms, we will refer the reader to the paper of Moret and Nualart
(2000).

By the isometry, and using Propositions \ref{corol2MN} and
\ref{lema10MN}
\begin{eqnarray*}
I_2&=&E \left(\sum_{t_i\in D_n,t_i\leq
t}f_{\Delta}^2(X_{t_{i}},t_{i})\int^{t_{i+1}}_{t_i}u_s^2ds\right)\\
&=&\sum_{t_i\in D_n,t_i\leq
t}E\left(f_{\Delta}^2(X_{t_{i}},t_{i})\int^{t_{i+1}}_{t_i}u_s^2ds\right)\\
&\leq&\sum_{t_i\in D_n,t_i\leq t}\int_{\R}f_{\Delta}^2(x,t_{i})dx
E\Big\vert \delta \left( \frac{ (\int_{t_i}^{t_{i+1}} u_s^2 ds )
DX_{t_i}}{ \Vert
DX_{t_i} \Vert_H^2} \right) \Big\vert\\
&\leq&C\sum_{t_i\in D_n,t_i\leq
t}\int_{\R}f_{\Delta}^2(x,t_{i})dxt_i^{-\frac12} \left( \sqrt{ E
\big\vert \int_{t_i}^{t_{i+1}} u_s^2 ds \big\vert^2} + \sqrt{ E
\int_0^{t_i} \left( D_t \left( \int_{t_i}^{t_{i+1}} u_s^2 ds
\right) \right)^2 dt} \right)\\
&=&C\sum_{t_i\in D_n,t_i\leq
t}\int_{\R}\sum_{k=1}^{m_1}\sum_{l=1}^{m_2}
f^2_{kl}I_{(x_k,x_{k+1}]}(x)I_{(s_l,s_{l+1}]}(t_i)
t_i^{-\frac12}(t_{i+1}-t_{i})dx\\
&=&C\int_{\R}\sum_{k=1}^{m_1}\sum_{l=1}^{m_2}
f^2_{kl}I_{(x_k,x_{k+1}]}(x)\left(\sum_{t_i\in D_n,t_i\leq
t}I_{(s_l,s_{l+1}]}(t_i)t_i^{-\frac12}(t_{i+1}-t_{i})\right)dx\\
&=&C\int_{\R}\sum_{k=1}^{m_1}\sum_{l=1}^{m_2}
f^2_{kl}I_{(x_k,x_{k+1}]}(x)\left(\int_0^t\sum_{t_i\in D_n,t_i\leq
t}I_{(s_l,s_{l+1}]}(t_i)t_i^{-\frac12}I_{(t_{i},t_{i+1}]}(s)ds\right)dx.
\end{eqnarray*}
Using the condition {\bf(M)} over the partitions, we have that, by
bounded convergence,
$$
\lim_{n\to\infty}\int_0^t\sum_{t_i\in D_n,t_i\leq t}
I_{(s_l,s_{l+1}]}(t_i)t_i^{-\frac12}I_{(t_{i},t_{i+1}]}(s)ds=
\int_0^t I_{(s_l,s_{l+1}]}(s)s^{-\frac12}ds,
$$
and then,
\begin{eqnarray}
\liminf_{n\to\infty}I_2&\leq&
C\int_{\R}\sum_{k=1}^{m_1}\sum_{l=1}^{m_2}
f^2_{kl}I_{(x_k,x_{k+1}]}(x)\int_0^t
I_{(s_l,s_{l+1}]}(s)s^{-\frac12}dsdx\nonumber\\
&=&C\int_0^t\int_{\R}f^2_{\Delta}(x,s)\frac1{\sqrt{s}}dxds\nonumber\\
&\leq&C\|f_{\Delta}\|^2.\label{aaa2}
\end{eqnarray}

On the other hand,
\begin{eqnarray}
I_1&=&E \left( \vert \sum_{t_i\in D_n,t_i\leq
t}f_{\Delta}(X_{t_{i+1}},t_{i+1})\left(X_{t_{i+1}}-X_{t_i}\right) \vert^2 \right)\nonumber\\
&=&E \left( \sum_{t_i\in D_n,t_i\leq
t}f_{\Delta}^2(X_{t_{i+1}},t_{i+1})\left(X_{t_{i+1}}-X_{t_i}\right)^2 \right) \nonumber\\
&&+2 E \left( \sum_{t_i,t_j\in D_n,t_i<t_j\leq
t}f_{\Delta}(X_{t_{i+1}},t_{i+1})f_{\Delta}(X_{t_{j+1}},t_{j+1})\left(X_{t_{i+1}}-X_{t_i}\right)
\left(X_{t_{j+1}}-X_{t_j}\right) \right) \nonumber\\
&:=&I_{1,1}+2I_{1,2}.\label{aaa3}
\end{eqnarray}

Following now the methods of Proposition 14 of Moret and Nualart
(1999) and using again  Propositions \ref{corol2MN} and
\ref{lema10MN} as we did in the study of $I_2$, we get that
$$I_{1,1}\leq C\sum_{t_i\in D_n,t_i\leq
t}\int_{\R}f_{\Delta}^2(x,t_{i})dxt_{i+1}^{-\frac12}(t_{i+1}-t_{i}).$$
By similar computations to those of the term $I_2$ we obtain that
\begin{equation}\liminf_{n\to\infty}I_{1,1}\leq
C\|f_{\Delta}\|^2.\label{aaa4}\end{equation}

Let us study now $I_{1,2}$. Using Lemma \ref{lemafita}, notice
that
\begin{eqnarray*}
I_{1,2}&=&E \left( \sum_{t_i,t_j\in D_n,t_i<t_j\leq
t}f_{\Delta}(X_{t_{i+1}},t_{i+1})f_{\Delta}(X_{t_{j+1}},t_{j+1})\left(X_{t_{i+1}}-X_{t_i}\right)
\left(X_{t_{j+1}}-X_{t_j}\right) \right) \\
&\leq&\sum_{t_i,t_j\in D_n,t_i<t_j\leq
t}E\left(f_{\Delta}(X_{t_{i+1}},t_{i+1})f_{\Delta}(X_{t_{j+1}},t_{j+1})C_{ij}\right).
\end{eqnarray*}
Following again the methods of the proof of Proposition 14 of
Moret and Nualart (1999) -more precisely, the proof of
inequalities (5.36) and (5.37)- the last expression is bounded by
\begin{eqnarray}
&&\sum_{t_i,t_j\in D_n,t_i<t_j\leq
t}E\left(f_{\Delta}^2(X_{t_{i+1}},t_{i+1})f_{\Delta}^2(X_{t_{j+1}},t_{j+1})\right)^{\frac12}E\left(C^2_{ij}\right)^{\frac12}\nonumber\\
&\leq&C\sum_{t_i,t_j\in D_n,t_i<t_j\leq
t}\frac{(t_{i+1}-t_i)(t_{j+1}-t_j)}{(t_{i+1}(t_{j+1}-t_{i+1}))^{\frac34}}\left(\int_{\R}f^2_{\Delta}(x,t_{i+1})dx\right)^{\frac12}
\left(\int_{\R}f^2_{\Delta}(x,t_{j+1})dx\right)^{\frac12}\nonumber\\
&\leq &C\sum_{t_i,t_j\in D_n,t_i<t_j\leq t}
\frac{(t_{i+1}-t_i)(t_{j+1}-t_j)}{(t_{i+1}(t_{j+1}-t_{i+1}))^{\frac{3}{4}}}
\left(\int_{\R}f^2_{\Delta}(x,t_{i+1})dx+
\int_{\R}f^2_{\Delta}(x,t_{j+1})dx\right)\nonumber\\
&:=&C(I_{1,2,1}+I_{1,2,2}).\label{aaa5}
\end{eqnarray}

Since
$$\sum_{t_j\in D_n,t_i<t_j\leq t}\frac{(t_{j+1}-t_j)}{(t_{j+1}-t_{i+1})^{\frac34}}=\int_{t_{i+1}}^t
\sum_{t_j\in D_n,t_i<t_j\leq
t}\frac{1}{(t_{j+1}-t_{i+1})^{\frac34}}I_{(t_j,t_{j+1}]}(s)ds,$$
we get that
\begin{eqnarray*}
I_{1,2,1}&\leq&\sum_{t_i\in D_n,t_i<t_j\leq t}
\frac{(t_{i+1}-t_i)}{t_{i+1}^{\frac{3}{4}}} \int_{t_{i+1}}^t
\frac{1}{(s-t_{i+1})^{\frac34}}ds
\left(\int_{\R}f^2_{\Delta}(x,t_{i+1})dx\right)\\
&\leq&C\sum_{t_i\in D_n,t_i<t_j\leq t}
\frac{(t_{i+1}-t_i)}{t_{i+1}^{\frac{3}{4}}}
\int_{\R}f^2_{\Delta}(x,t_{i+1})dx.
\end{eqnarray*}
And this clearly yields that
\begin{equation}\liminf_{n\to\infty}I_{1,2,1}\leq
C\int_0^t\int_{\R}f^2_{\Delta}(x,s)\frac1{s^{\frac34}}ds
=C\|f_{\Delta}\|^2.\label{aaa6}\end{equation}

Finally we have to consider $I_{1,2,2}$. First of all, notice that
\begin{eqnarray*}
I_{1,2,2}&:=&\sum_{t_i,t_j\in D_n,t_i<t_j\leq t}
\frac{(t_{i+1}-t_i)(t_{j+1}-t_j)}{(t_{i+1}(t_{j+1}-t_{i+1}))^{\frac{3}{4}}}
\left(\int_{\R}f^2_{\Delta}(x,t_{j+1})dx\right)\\
&=&\sum_{t_i,t_j\in D_n,t_i<t_j\leq t}
\frac{(t_{i+1}-t_i)(t_{j+1}-t_j)}{(t_{i+1}(t_{j+1}-t_{i+1}))^{\frac{3}{4}}}
\left(\int_{\R} \sum_{k=1}^{m_1}\sum_{l=1}^{m_2}
f^2_{kl}I_{(x_k,x_{k+1}]}(x)I_{(s_l,s_{l+1}]}(t_{j+1})
dx\right)\\
&\leq&\int_{\R}
\sum_{k=1}^{m_1}\sum_{l=1}^{m_2}f^2_{kl}I_{(x_k,x_{k+1}]}(x)\sum_{t_i\in
D_n,t_i< t} \frac{(t_{i+1}-t_i)}{t_{i+1}^{\frac{3}{4}}}
\\&&\times\sum_{t_j\in D_n,t_i<t_j\leq
t}\int_{t_{i+1}}^t\frac1{(s-t_{i+1})^{\frac34}}
I_{(s_l,s_{l+1}]}(t_{j+1})I_{(t_{j},t_{j+1}]}(s)dsdx
\end{eqnarray*}
From the obvious inequality
$$I_{(s_l,s_{l+1}]}(t_{j+1})I_{(t_{j},t_{j+1}]}(s)\leq I_{(s_l,s_{l+1}]}(s)I_{(t_{j},t_{j+1}]}(s)
+I_{(t_j,s_l]}(s)I_{(t_{j},t_{j+1}]}(s_l),$$
we obtain the bound
\begin{equation}I_{1,2,2}\leq I_{1,2,2,1}+I_{1,2,2,2},\label{aaa9}\end{equation}
where
\begin{eqnarray*}
I_{1,2,2,1}&=&\int_{\R}
\sum_{k=1}^{m_1}\sum_{l=1}^{m_2}f^2_{kl}I_{(x_k,x_{k+1}]}(x) \\
& & \qquad \times \sum_{t_i\in D_n,t_i< t}
\frac{(t_{i+1}-t_i)}{t_{i+1}^{\frac{3}{4}}} \sum_{t_j\in
D_n,t_i<t_j\leq t}\int_{t_{i+1}}^t\frac1{(s-t_{i+1})^{\frac34}}
I_{(s_l,s_{l+1}]}(s)I_{(t_{j},t_{j+1}]}(s)dsdx\\
I_{1,2,2,2}&=&\int_{\R}
\sum_{k=1}^{m_1}\sum_{l=1}^{m_2}f^2_{kl}I_{(x_k,x_{k+1}]}(x)\\
& & \qquad \times \sum_{t_i\in D_n,t_i< t}
\frac{(t_{i+1}-t_i)}{t_{i+1}^{\frac{3}{4}}} \sum_{t_j\in
D_n,t_i<t_j\leq t}\int_{t_{i+1}}^t\frac1{(s-t_{i+1})^{\frac34}}
I_{(t_j,s_{l}]}(s)I_{(t_{j},t_{j+1}]}(s_l)dsdx.
\end{eqnarray*}

Now, since we can write
$$
I_{1,2,2,1} = \int_{\R} \sum_{t_i\in D_n,t_i< t}
\frac{(t_{i+1}-t_i)}{t_{i+1}^{\frac{3}{4}}}
\int_{t_{i+1}}^tf^2_{\Delta}(x,s)\frac1{(s-t_{i+1})^{\frac34}}
dsdx,
$$
using an argument of bounded convergence we have that
\begin{eqnarray}
\liminf_{n\to\infty}I_{1,2,2,1}& \leq &\int_{\R}\int_0^t
\frac{1}{u^{\frac{3}{4}}}
\int_{u}^tf^2_{\Delta}(x,s)\frac1{(s-u)^{\frac34}} dsdudx\nonumber\\
&=&\int_{\R}\int_0^t
f^2_{\Delta}(x,s)\int_{0}^s\frac{1}{u^{\frac{3}{4}}}
\frac1{(s-u)^{\frac34}} dudsdx\nonumber\\
&\leq&C\int_{\R}\int_0^t f^2_{\Delta}(x,s)\frac{1}{s^{\frac12}}
dsdx\nonumber\\
&\leq&C\|f_{\Delta}\|^2\label{aaa10}.
\end{eqnarray}

On the other hand, observe that fixed $l$, there exists only one
$j$ (that we will denote by $j(l)$) such that $t_{j(l)}<s_l\leq
t_{j(l)+1}$. So,
\begin{eqnarray*}
&&\sum_{t_i,t_j\in D_n,t_i<t_j\leq t}
\frac{(t_{i+1}-t_i)}{t_{i+1}^{\frac{3}{4}}}
\int_{t_{i+1}}^t\frac1{(s-t_{i+1})^{\frac34}}
I_{(t_{j},s_{l}]}(s)I_{(t_{j},t_{j+1}]}(s_l)ds\\
&\leq&\sum_{t_i\in D_n,t_i<t_{j(l)}\leq t}
\frac{(t_{i+1}-t_i)}{t_{i+1}^{\frac{3}{4}}}
\int_{t_{j(l)}}^{t_{j(l)+1}}\frac1{(s-t_{i+1})^{\frac34}}ds.
\end{eqnarray*}
Now, using that for $i<j(l)$
$$\int_{t_{j(l)}}^{t_{j(l)+1}}\frac1{(s-t_{i+1})^{\frac34}}ds
\leq \int_{t_{j(l)}}^{t_{j(l)+1}}\frac1{(s-t_{j(l)})^{\frac34}}ds
\leq4|D_n|^{\frac14},$$ and that
$$\sum_{t_i\in D_n,t_i<t_{j(l)}\leq t}
\frac{(t_{i+1}-t_i)}{t_{i+1}^{\frac{3}{4}}}\leq\int_0^1\frac1{s^{\frac34}}ds<\infty,$$
we obtain easily that
\begin{equation}\lim_{n\to\infty}I_{1,2,2,2}=0.\label{aaa7}
\end{equation}

So, putting together (\ref{aaa1})-(\ref{aaa7}), we have proved
(\ref{simpli}).

 Now, given $f\in\mathcal
H$, let us consider  $\{f_n\}_{n\in\N}$ a sequence of elementary
functions converging to  $f$ in $\mathcal H$, and we define
$$\int_0^t\int_{\R}f(x,s)dL_s^x=L^1-\lim_{n\to\infty}\left(\int_0^t\int_{\R}f_{n}(x,s)dL_s^x\right).$$
Clearly, this limit exists. Indeed, for any $\varepsilon>0$ there
exists $n_0$ such that for any $n,m\geq n_0$,
$\|f_n-f_m\|<\varepsilon$ and using inequality (\ref{simpli}) we
obtain that
\begin{eqnarray*}
E\left|\int_0^t\int_{\R}f_n(x,s)dL_s^x-\int_0^t\int_{\R}f_{m}(x,s)dL_s^x\right|
=E\left|\int_0^t\int_{\R}(f_n(x,s)-f_m(x,s))dL_s^x\right|\leq\|f_n-f_m\|<\varepsilon.
\end{eqnarray*}
Moreover, using again inequality (\ref{simpli}), it is clear that the
definition does not depend on the choice of the sequence $(f_n)$.
Indeed, given $(f_n^1)_{n\in\N}$ and $(f_n^2)_{n\in\N}$ two
sequences converging to $f$ in $\mathcal H$, we have
\begin{eqnarray*}
E\left(\left|\int_0^t\int_{\R}f_n^1(x,s)dL_s^x-\int_0^t\int_{\R}f_{n}^2(x,s)dL_s^x\right|\right)\leq\|f_n^1-f_n^2\|\leq\|f_n^1-f\|+\|f-f_n^2\|,
\end{eqnarray*}
that goes to zero when $n$ tends to infinity.

\hfill
$\Box$
\end{dem}

\begin{obs}\label{obsC}
If $f$ satisfies condition {\bf (C)}, from Theorem \ref{teoM} we
know that the quadratic covariation $\left[f(X,\cdot),X\right]$
exists. Moreover, if $f\in\mathcal H$, from  the uniqueness of the
extension in the construction of the integral in  Theorem
\ref{teo2} we get that
$$\int_0^t\int_{\R}f(x,s)dL_s^x=-\left[f(X,\cdot),X\right]_t.$$
\end{obs}

The following results is an obvious consequence of Theorem
\ref{teoM} and Remark \ref{obsC}.

\begin{cor}\label{elcor}
Let $u$ be a process satisfying {\bf (H1)} and  {\bf (H2)}. Set
$X=\int_0^t us dW_s$. Consider a sequence of partitions $D_n$ of
partitions of $[0,1]$ verifying conditions {\bf (M)}. Let $F(x,t)$
be an absolutely continuous function in $x$ such that the partial
derivative $f(\cdot,t)$ satisfies {\bf (C)}. Then, if
$f\in\mathcal H$, we have the following extension for the
It{\^o}'s formula:
$$F(X_t,t)=F(0,0)+\int_0^tf(X_s,s)dX_s-\frac12\int_0^t\int_{\R}f(x,s)dL_s^x+\int_0^tF(X_s,ds).$$
\end{cor}

\section{It{\^o}'s formula extension}

Now we can state the main result of this paper.

\begin{teo}\label{teo3} $ $

Hypotheses over the martingale:
\begin{enumerate}\item Let $u$ be an adapted process satisfying {\bf (H1)} and  {\bf (H2)}. Set
$X=\int_0^t u_s dW_s$.
\end{enumerate}

Hypotheses over the function:

 \begin{enumerate} \item Let
$F$ be a function defined on $\R\times[0,1]$ such that $F$ admits
first order Radon-Nikodym derivatives with respect to each
parameter. \item Assume that these derivatives satisfy that for
every $A\in\R$,
\begin{eqnarray*}
\int_0^1\int_{-A}^{A}\left|\frac{\partial F}{\partial
t}(x,s)\right|dx \frac1{\sqrt{s}} ds&<&+\infty \\
\int_0^1\int_{-A}^{A}\left(\frac{\partial F}{\partial
x}(x,s)\right)^2dx \frac1{\sqrt{s}} ds&<&+\infty.
\end{eqnarray*}
\end{enumerate}
Then, for all $t\in[0,1]$,
\begin{eqnarray*}
F(X_t,t)=F(0,0)+\int_0^t\frac{\partial F}{\partial
x}(X_s,s)dX_s+\int_0^t\frac{\partial F}{\partial
t}(X_s,s)ds-\frac12\int_0^t\int_{\R}\frac{\partial F}{\partial
x}(x,s)dL_s^x.
\end{eqnarray*}
\end{teo}

\begin{dem}
Using localization arguments we can assume that $F$ has compact
support and
\begin{eqnarray*}
\int_0^1\int_{\R}\left|\frac{\partial F}{\partial
t}(x,s)\right|dx \frac1{\sqrt{s}} ds&<&+\infty \\
\int_0^1\int_{\R}\left(\frac{\partial F}{\partial
x}(x,s)\right)^2dx \frac1{\sqrt{s}} ds&<&+\infty.
\end{eqnarray*}

\smallskip

Let $g\in\mathcal C^{\infty}$ be a function with compact support
from $\R$ to $\R^{+}$ such that  $\int_{\R}g(s)ds=1$. We define,
for any $n\in\N$,
$$g_n(s)=ng(ns)$$
and
$$F_n(x,t)=\int_0^1\int_{\R}F(y,s)g_n(t-s)g_n(x-y)dyds.$$
Then  $F_n\in\mathcal C^{\infty}(\R\times[0,1])$. Hence, by the
usual It{\^o}'s formula, for every $\varepsilon>0$, we can write
\begin{eqnarray}\label{nova}
F_n(X_t,t)=F_n(X_{\varepsilon},\varepsilon)+\int_{\varepsilon}^t\frac{\partial
F_n}{\partial x}(X_s,s)dX_s+\int_{\varepsilon}^t\frac{\partial
F_n}{\partial t}(X_s,s)ds+\frac12\int_{\varepsilon}^t
u_s^2\frac{\partial^2 F_n}{\partial x^2}(X_s,s)ds.
\end{eqnarray}

\smallskip

Using the arguments of Az{\'e}ma {\it et al.} (1998) we will study
the convergence of (\ref{nova}).

\smallskip

 Since $F$ is
a continuous function with compact support,   it is easy to check
that $(F_n(X_t,t))_{n\in\N}$ converges in probability to
$F(X_t,t)$.

\smallskip

On the other hand
$$\int_0^1\int_{\R} \left| \frac{\partial
F}{\partial t}(x,s) \right| dxds\leq\int_0^1\int_{\R}\left|\frac{\partial
F}{\partial t}(x,s)\right|dx \frac1{\sqrt{s}} ds<+\infty.$$ Hence,
$\frac{\partial F}{\partial t}\in L^1(\R\times[0,1])$. Under our
hypothesis over the martingale $X$, it follows from Proposition
\ref{2.6MN} and Lemma \ref{lema10MN} that for any $t \in [0,1]$,
the random variable $X_t$ is absolutely continuous with density
$p_t$ satisfying the estimate
$$p_t(x)\leq\frac{C}{\sqrt{t}}.$$
Then, it is easy to see that
$\left(\int_{\varepsilon}^t\frac{\partial F_n}{\partial
t}(X_s,s)ds\right)_{n\in\N}$ converges in probability to
$\left(\int_{\varepsilon}^t\frac{\partial F}{\partial
t}(X_s,s)ds\right)$. Indeed,
\begin{eqnarray*}
E\left(\left|\int_{\varepsilon}^t\left(\frac{\partial
F_n}{\partial t}(X_s,s)-\frac{\partial F}{\partial
t}(X_s,s)\right)ds\right|\right)
&\leq&\int_{\varepsilon}^t\int_{\R}\left|\frac{\partial
F_n}{\partial t}(x,s)-\frac{\partial F}{\partial t}(x,s)
\right|p_s(x)dxds\\
&\leq&C\int_{\varepsilon}^t\int_{\R}\left|\frac{\partial
F_n}{\partial t}(x,s)-\frac{\partial F}{\partial t}(x,s)
\right|\frac1{\sqrt{s}}dxds\\
&\leq&\frac{C}{\sqrt{\varepsilon}}\int_{\varepsilon}^t\int_{\R}\left|\frac{\partial
F_n}{\partial t}(x,s)-\frac{\partial F}{\partial t}(x,s)
\right|dxds,
\end{eqnarray*}
that goes to zero, when $n$ tends to infinity, since
$\frac{\partial F}{\partial t}\in L^1(\R\times[0,1])$ and
$$\frac{\partial F_n}{\partial t}(x,t)=\int_0^1\int_{\R}\frac{\partial F}{\partial t} (y,s)g_n(t-s)g_n(x-y)dyds.$$

\smallskip

Similarly, we can prove that
$\left(\int_{\varepsilon}^t\frac{\partial F_n}{\partial
x}(X_s,s)dX_s\right)_{n\in\N}$ converges in probability to \break
$\left(\int_{\varepsilon}^t\frac{\partial F}{\partial
x}(X_s,s)dX_s\right)$. Indeed, using the same arguments we get
that $\frac{\partial F}{\partial x}\in L^2(\R\times[0,1]).$ Then,
\begin{eqnarray*}
&&E\left(\left|\int_{\varepsilon}^t\left(\frac{\partial
F_n}{\partial x}(X_s,s)-\frac{\partial F}{\partial
x}(X_s,s)\right)dX_s\right|^2\right)
\\&&\quad = E\left(\left|\int_{\varepsilon}^t\left(\frac{\partial
F_n}{\partial x}(X_s,s)-\frac{\partial F}{\partial
x}(X_s,s)\right)u_sdW_s\right|^2\right)\\
&&\quad = E\left(\int_{\varepsilon}^t\left(\frac{\partial
F_n}{\partial x}(X_s,s)-\frac{\partial F}{\partial
x}(X_s,s)\right)^2u_s^2ds\right).
\end{eqnarray*}
Following the same ideas of Proposition 12 in Moret and Nualart
(2000),  Proposition \ref{corol2MN} and Lemma \ref{lema10MN} yield
the following  bound for the last expression
\begin{eqnarray}
&&C \int_{\varepsilon}^t\int_{\R}\left(\frac{\partial
F_n}{\partial x}(x,s)-\frac{\partial F}{\partial x}(x,s)
\right)^2\frac1{\sqrt{s}}dxds\nonumber\\
&\leq&\frac{C}{\sqrt{\varepsilon}}
\int_{\varepsilon}^t\int_{\R}\left(\frac{\partial F_n}{\partial
x}(x,s)-\frac{\partial F}{\partial x}(x,s)
\right)^2dxds\label{ladarrera}
\end{eqnarray}
that goes to zero when $n$ tends to infinity, since
$\frac{\partial F}{\partial x}\in L^2(\R\times[0,1])$ and
$$\frac{\partial F_n}{\partial x}(x,t)=\int_0^1\int_{\R}\frac{\partial F}{\partial x} (y,s)g_n(t-s)g_n(x-y)dyds.$$

\smallskip

So, letting $n$ to infinity in (\ref{nova}), we get that the
sequence
$$\left(\frac12\int_{\varepsilon}^tu_s^2\frac{\partial^2
F_n}{\partial x^2}(X_s,s)ds\right)_{n\in\N}$$ converges in
probability to
$$F(X_t,t)-F(X_{\varepsilon},\varepsilon)-\int_{\varepsilon}^t\frac{\partial
F}{\partial x}(X_s,s)dX_s-\int_{\varepsilon}^t\frac{\partial
F}{\partial t}(X_s,s)ds.$$ But, since $\frac{\partial
F_n}{\partial x}(x,s)I_{(\varepsilon,t)}(s)\in\mathcal H$, from
Theorem \ref{teoM} and Corollary \ref{elcor}, we get that
\begin{eqnarray*}
\int_{\varepsilon}^tu_s^2\frac{\partial^2 F_n}{\partial
x^2}(X_s,s)ds&=&\left[\frac{\partial F_n}{\partial
x}(X,\cdot),X\right]_t-\left[\frac{\partial F_n}{\partial
x}(X,\cdot),X\right]_{\varepsilon}\\
&=&-\int_0^1\int_{\R}\frac{\partial F_n}{\partial
x}(x,s)I_{(\varepsilon,t)}(s)dL_s^x.
\end{eqnarray*}

\smallskip

The next step of the proof is to check that $\left( \frac{\partial
F_n}{\partial
x}(x,s)I_{(\varepsilon,t)}(s),x\in\R,s\in[0,1]\right)_{n\in\N}$
converges in $\mathcal H$ to $\big( \frac{\partial F}{\partial
x}(x,s)I_{(\varepsilon,t)}(s),x\in\R, s\in[0,1]\big)$. It suffices
to notice that,
$$\int_{\varepsilon}^t\int_{\R}\left(\frac{\partial F_n}{\partial
x}(x,s)-\frac{\partial F}{\partial
x}(x,s)\right)^2\frac1{s^{\frac34}}dxds\leq
\frac1{{\varepsilon}^{\frac34}}\int_{\varepsilon}^t\int_{\R}\left(\frac{\partial
F_n}{\partial x}(x,s)-\frac{\partial F}{\partial
x}(x,s)\right)^2dxds$$ that converges to zero when $n$ tends to
infinity. Then, we clearly have proved that
$$\left(\int_0^1\int_{\R}\frac{\partial F_n}{\partial
x}(x,s)I_{(\varepsilon,t)}(s)dL_s^x\right)_{n\in\N}$$ converges in
$L^1$ to $\int_0^1\int_{\R}\frac{\partial F}{\partial
x}(x,s)I_{(\varepsilon,t)}(s)dL_s^x$.

\smallskip

So, we have that for any $\varepsilon>0$
\begin{eqnarray}\label{eq4}
F(X_t,t)=F(X_{\varepsilon},\varepsilon)+\int_{\varepsilon}^t\frac{\partial
F}{\partial x}(X_s,s)dX_s+\int_{\varepsilon}^t\frac{\partial
F}{\partial t}(X_s,s)ds-\frac12\int_0^1\int_{\R}\frac{\partial
F}{\partial x}(x,s)I_{(\varepsilon,t)}(s)dL_s^x.
\end{eqnarray}
The last steep is to let $\varepsilon$ to zero. But  we need to
check that the limit of the stochastic integral exists. Actually,
it is enough to show that
$$E\left(\left|\int_0^t\frac{\partial F}{\partial
t}(X_s,s)ds\right|\right)<\infty$$ and that
$$E\left(\int_0^t\frac{\partial F}{\partial x}(X_s,s)dX_s\right)^2<\infty.$$
But,
$$E\left(\left|\int_0^t\frac{\partial F}{\partial
t}(X_s,s)ds\right|\right)\leq C
\int_0^1\int_{\R}\left|\frac{\partial F}{\partial t}(x,s)\right|dx
\frac1{\sqrt{s}} ds<+\infty.$$

On the other hand, following the same type of arguments that in
(\ref{ladarrera}), we are able to write
\begin{eqnarray*}
&&E\left(\int_0^t\frac{\partial F}{\partial x}(X_s,s)dX_s\right)^2
= E \left(\int_0^t\left(\frac{\partial F}{\partial
x}(X_s,s)\right)^2u_s^2ds\right)\\ \quad & & \leq
C \int_\R \int_0^t\left(\frac{\partial F}{\partial
x}(x,s)\right)^2 \frac1{\sqrt{s}}dsdx <\infty.
\end{eqnarray*}

Letting $\varepsilon$ to zero, the proof is finished.

\hfill $\Box$
\end{dem}

\begin{obs}\label{obs6}
Notice that under the hypotheses   of Theorem \ref{teo3}, it is
possible that $\frac{\partial F}{\partial x}$  does not belong to
the space $\mathcal H$. In this case, using the localization
arguments, we can always assume that $\left( \frac{\partial
F}{\partial x}(x,s)I_{(\varepsilon,t)}(s),\, x\in\R,
s\in[0,1]\right)$ belongs to $\mathcal H$ for any $\varepsilon>0$
and we can define
$$\int_0^t\int_{\R}\frac{\partial
F}{\partial x}(x,s)dL_s^x=\lim_{\varepsilon\to
0}\int_0^1\int_{\R}\frac{\partial F}{\partial
x}(x,s)I_{(\varepsilon,t)}(s)dL_s^x.$$ This limit exists in
probability since all the other limits in (\ref{eq4}) exist.
\end{obs}

\section*{Acknowledgements}
This work was partially supported by DGES Grants MTM2006-01351 (Carles Rovira) and
MTM2006-06427 (Xavier Bardina).

\begin{description}

\item[] {Az{\'e}ma, J., Jeulin, T., Knight, F., Yor, M. (1998)}
{Quelques calculs de compensateurs impliquant l'injectivit{\'e} de
certains processus croissants.} {\it S{\'e}minaire de Probabilit{\'e}s,
XXXII, Lecture Notes in Math., {\bf 1686},} {316-327.}

\item[] {Bardina, X., Jolis, M. (1997)} {An extension of It\^{o}'s
formula for elliptic martingale processes.} {\it Stochastic
Process. Appl.} {\bf 69} {(1), 83-109.}

\item[] {Bardina, X., Jolis, M. (1997)} {Estimations of the density of hypoelliptic diffusion processes with applications to an extension of It\^{o}'s
formula.} {\it J. Theoret. Probab.} {\bf 15} {(1), 223-247.}

\item[] {Bardina, X., Rovira, C. (2007)} {On It\^{o}'s
formula for elliptic diffusion processes.} {\it Bernoulli} {\bf
13} {(3), 820-830.}

\item[] {Bouleau, N., Yor, M. (1981)} {Sur la variation
quadratique des temps locaux de certaines semimartingales.} {\it
C. R. Acad. Sci. Paris S{\'e}r. I Math.} {\bf 292} {(9), 491-494.}

\item[] {Di Nunno, G., Meyer-Brandis, T., {\O}ksendal, B., Proske, F. (2005)} Malliavin calculus and anticipative It\^o formulae for L\'evy processes.  {\it Infin. Dimens. Anal. Quantum Probab. Relat. Top.} {\bf  8}  (2), 235--258.

\item[] {Dupoiron, K., Mathieu, P., San Martin, J.  (2004)}
{Formule d'It{\^o} pour des martingales uniform{\'e}ment
elliptiques, et processus de Dirichlet.} {\it Potential Anal.}
{\bf 21} {(1), 7-33.}

\item[] {Eisenbaum, N. (2000)} {Integration with respect to local
time.} {\it Potential Anal.} {\bf 13} {(4), 303-328.}

\item[] {Eisenbaum, N.  (2001)} {On It{\^o}'s formula of F{\"o}llmer and
Protter.} {\it S{\'e}minaire de Probabilit{\'e}s, XXXV, Lecture Notes in
Math. {\bf 1755},} {390-395.}

\item[]  Flandoli, F., Russo, F., Wolf, J. (2004)   Some SDEs with distributional drift. II. Lyons-Zheng structure, It\^{o}'s
formula and semimartingale characterization.  {\it Random Oper. Stochastic Equations} {\bf  12} (2), 145--184.

\item[] {F{\"o}llmer, H., Protter, P., Shiryayev, A.N. (1995)}
{Quadratic covariation and an extension of It{\^o}'s formula.} {\it
Bernoulli} {\bf 1} {(1-2), 149-169.}

\item[] {Ghomrasni, R., Peskir, G. (2003)} {Local time-space
calculus and extensions of It{\^o}'s formula.} {\it High
dimensional probability, III (Sandjberg, 2002), Progr. Probab.,}
{\bf 55}, {177-192.}



\item[] {Moret, S. (1999)} {Ph.D. thesis: ``Generalitzacions de la Formula d'It\^o i estimacions per martingales".} {\it Universitat de Barcelona}.

\item[] {Moret, S., Nualart, D. (2000)} {Quadratic Covariation and It\^o's Formula for Smooth
Nondegenerate Martingales.} {\it Journal of Theoretical Probability} {\bf 13}, {193-224.}

\item[] {Nualart, D. (2006)}  {\it Malliavin Calculus and Related Topics. Second edition.} Springer-Verlag.

\item[] {Nualart, D. (1998)}  Analysis on Wiener space and anticipating stochastic calculus and Related Topics.
In: \'Ecole d'\'et\'e de Saint-Flour XXV. {\it Lect. Notes in Math.}
{\bf 1690}, {123-227.}



\end{description}

\end{document}